\documentclass[journal]{IEEEtran}

\usepackage{amsmath,amssymb,amsthm}
\usepackage{graphicx}
\usepackage{cite}
\usepackage{url}
\usepackage{booktabs}
\usepackage{algorithm}
\usepackage{algpseudocode}
\usepackage{xcolor}

\graphicspath{{figures/}}

\newtheorem{theorem}{Theorem}
\newtheorem{proposition}{Proposition}

\newcommand{\bx}{\boldsymbol{x}}
\newcommand{\bX}{\boldsymbol{X}}
\newcommand{\bW}{\boldsymbol{W}}
\newcommand{\bQ}{\boldsymbol{Q}}
\newcommand{\bc}{\boldsymbol{c}}
\newcommand{\bb}{\boldsymbol{b}}
\newcommand{\bu}{\boldsymbol{u}}

\newcommand{\bt}{\boldsymbol{t}}
\newcommand{\bv}{\boldsymbol{v}}
\newcommand{\be}{\boldsymbol{e}}
\newcommand{\bz}{\boldsymbol{z}}
\newcommand{\bq}{\boldsymbol{q}}
\newcommand{\bY}{\boldsymbol{Y}}
\newcommand{\bB}{\boldsymbol{B}}
\newcommand{\bF}{\boldsymbol{F}}
\newcommand{\bI}{\boldsymbol{I}}
\newcommand{\bU}{\boldsymbol{U}}
\newcommand{\bV}{\boldsymbol{V}}
\newcommand{\bD}{\boldsymbol{D}}
\newcommand{\bPhi}{\boldsymbol{\Phi}}
\newcommand{\blam}{\boldsymbol{\lambda}}

\newcommand{\B}{\mathcal{B}}
\newcommand{\Hb}{\mathbb{H}}
\newcommand{\C}{\mathbb{C}}
\newcommand{\R}{\mathbb{R}}
\newcommand{\tr}{\operatorname{tr}}
\newcommand{\sgn}{\operatorname{sgn}}
\newcommand{\re}{\operatorname{Re}}
\newcommand{\lead}{\operatorname{lead}}

\begin{document}

\title{Radar Waveform Design via Convex\\ Semidefinite Relaxation}

\author{Samuel~Pinilla,~\IEEEmembership{Member,~IEEE}
\thanks{S. Pinilla is with Diamond Light Source, Harwell Science Campus, Chilton-Didcot
OX11 0DE, U.K. (e-mail: spinilla@ieee.org).}}

\markboth{IEEE Signal Processing Letters}
{Pinilla: Radar Waveform Design via Convex Semidefinite Relaxation}

\maketitle

\begin{abstract}
The ambiguity function (AF) of a radar waveform describes its ability to resolve targets in range and Doppler. Designing a waveform that attains a prescribed AF magnitude is a phase retrieval problem with four trivial ambiguities, and a time- or band-limited pulse is determined up to these ambiguities by at most three times as many AF samples as it has time samples or frequencies. Existing solvers for the narrowband AF, however, are non-convex: \textit{they cannot guarantee that they reach the global minimizer, which is unique only up to the trivial ambiguities, so the recovered pulse depends on the starting point}. This letter addresses this gap by lifting the problem to positive semidefinite matrices supported on the pulse band, where each AF sample becomes the magnitude of a linear functional. In the lifted space, recovery is formulated as a convex program, and a dual certificate provides a sufficient condition for exactness. The proposed algorithm, \textit{wavex}, combines an initialization by relaxed reflections, convex semidefinite surrogate steps, and a local polish. Numerical experiments show that \textit{wavex} accurately recovers band-limited pulses from complete AFs and from AFs with most of their samples removed, provided the remaining samples determine the pulse, and that its accuracy degrades gracefully as the noise increases. \textit{wavex} is built on the open-source library \textit{invexapi}, and the code that reproduces every result is available at \texttt{\url{https://github.com/samuelpinilla/waveX}}.
\end{abstract}


\section{Introduction}

\IEEEPARstart{T}{he} matched filter~\cite{north1963} is optimal for an echo of known shape but does not describe the response to any other echo. Woodward filled this gap~\cite{woodward1953} with the ambiguity function (AF), the filter output versus delay and Doppler mismatch, whose properties were developed in~\cite{siebert1956,price1965}. For a discrete signal $\bx\in\C^N$, the conventional AF is
\begin{equation}
  \boldsymbol{A}[\tau,f] = \left\lvert\sum_{n=0}^{N-1} \bx[n]\,\overline{\bx[n-\tau]}\,e^{-2\pi i nf/N}\right\rvert^2,
  \label{eq:afdiscrete}
\end{equation}
with $\tau,f=0,\dots,N-1$ and indices modulo $N$. Its bilinear form traces back to Wigner~\cite{wigner1932} and Ville~\cite{ville1948}, and the AF was later interpreted through the Heisenberg group~\cite{schempp1984,auslander1985} and extended to wideband~\cite{kelly1965,speiser1967} and multiple-input multiple-output~\cite{sanantonio2007} radar.

The mainlobe width of the AF sets the range and velocity resolution, and its sidelobes indicate how strongly a target in one delay-Doppler cell can mask a target in another, which makes the AF the standard criterion for matching a waveform to an application~\cite{wilcox1960}. A response concentrated at a single delay-Doppler point is impractical, since a target would rarely fall exactly on it~\cite{jankiraman2007,rihaczek1996}, so waveform design seeks a thumbtack AF, with a narrow central peak and low sidelobes. Finding the signal that produces a prescribed AF is the inverse of this task; in an early result, de~Buda~\cite{debuda1970} showed that an AF bounded by a Hermite function belongs to a Hermite-function signal. The following \emph{trivial ambiguities} of $\bx$ leave~\eqref{eq:afdiscrete} unchanged:
\begin{description}
  \item[T1] Rotated signal $e^{i\phi}\bx$ for some $\phi\in\R$.
  \item[T2] Translated signal $\bx[n-a]$ for some $a\in\R$.
  \item[T3] Reflected signal $\bx[-n]$.
  \item[T4] Modulated signal $e^{ibn}\bx[n]$ for some $b\in\R$.
\end{description}
Since the AF is phaseless, the inverse problem is a bivariate phase retrieval (PR) problem~\cite{shechtman2015,jaming2010}, related to short-time Fourier and frequency-resolved optical gating (FROG) PR~\cite{bendory2020frog,pinilla2019frog}. The conjugation in the AF adds a fourth trivial ambiguity to those of FROG~\cite{pinilla2024iaae}.

For almost all band- or time-limited pulses, at most three times as many AF samples as frequencies or time samples determine the pulse up to T1--T4~\cite{pinilla2024iaae,pinilla2021banraw}. Existing solvers are non-convex: the trust-region method of~\cite{pinilla2024iaae}, minimizing a smoothed amplitude loss~\cite{pinilla2018smoothing,pinilla2018ssp}, guarantees convergence only to a critical point. Lacking a certificate of global optimality, its output cannot be accepted as the unique solution~\cite{pinilla2024invex}, a concern across phaseless problems~\cite{pinilla2023bispectrum,pinilla2021asilomar,jerez2020detection,mishra2022group}. WaveMax~\cite{pinilla2021wavemax,pinilla2025wavemax}, a convex program in the line of~\cite{candes2013phaselift,guerrero2020guarantees,bacca2020superres,pinilla2023spm}, is posed on the AF of the fractional Fourier transform of the pulse, whose samples differ from~\eqref{eq:afdiscrete}, so its program and guarantees do not transfer directly. A convex formulation for the conventional AF, the one used by most radar systems~\cite{levanon2004}, is still missing.

This letter fills that gap with \textit{wavex}, a convex method based on lifting the AF. Its contributions are: (i) a band-restricted lifting in which every sample of~\eqref{eq:afdiscrete} is the magnitude of a linear functional of $\bX=\bx\bx^H$; (ii) a convex program over this lifting, with a dual certificate as a sufficient condition for exact recovery; (iii) a new initialization by relaxed averaged reflections from a single random start; (iv) convex surrogate steps that update the lifted matrix by proximal-gradient iterations, followed by a local polish; and (v) experiments on complete, decimated and noisy AFs. The code is public~\cite{wavex} and built on \textit{invexapi}~\cite{invexapi}.

\emph{Notation:} Boldface lowercase and uppercase letters denote vectors and matrices. $\Hb^N$ denotes the $N\times N$ Hermitian matrices, $\langle \bU,\bV\rangle = \tr(\bU^H\bV)$, $\|\cdot\|_F$, $\|\cdot\|_*$ and $\|\cdot\|_2$ the Frobenius, nuclear and spectral norms of a matrix, $\odot$ the entrywise product, $\sgn(z)=z/|z|$ (with $\sgn(0)=0$), and $\lead(\bW) = \sqrt{\lambda_1}\,\bv_1$ the scaled leading eigenvector of $\bW$.

\section{Problem Formulation}

Let $\bx\in\C^N$ be the pulse. A delay by $\tau\in\R$ is applied by the operator $\bPhi_\tau = \bF^H\operatorname{diag}\big(e^{-2\pi i l\tau/N}\big)\bF$, where $\bF$ is the unitary DFT and $l=-N/2,\dots,N/2-1$ indexes its frequencies. For integer $\tau$ it is the cyclic shift of~\eqref{eq:afdiscrete}, $(\bPhi_\tau\bx)[n] = \bx[n-\tau]$, and for real $\tau$ it defines the translation in T2. With the matrices $\bB_{\tau,f} = \bPhi_\tau^H\operatorname{diag}\big(e^{-2\pi i nf/N}\big)$, $n=0,\dots,N-1$, the sum in~\eqref{eq:afdiscrete} is the quadratic form $\bx^H\bB_{\tau,f}\bx$, so that
\begin{equation}
  \boldsymbol{A}[\tau,f] = \big|\bx^H\bB_{\tau,f}\,\bx\big|^2 = \big|\tr(\bB_{\tau,f}\bX)\big|^2,
  \quad (\tau,f)\in\Omega ,
  \label{eq:af}
\end{equation}
where $\bX=\bx\bx^H$ is the lifted pulse and $\Omega$ is the set of $m$ measured delay-Doppler pairs. Every AF sample is thus the squared magnitude of a linear functional of $\bX$.

Let $\B:\Hb^N\to\C^m$ stack the functionals $\tr(\bB_{\tau,f}\,\cdot\,)$ for $(\tau,f)\in\Omega$, and let $\bb\in\R^m$ stack the amplitudes $\bb[\tau,f]=\sqrt{\boldsymbol{A}[\tau,f]}$. Unlike classical lifting, where $|\boldsymbol{a}^H\bx|^2 = \tr(\boldsymbol{a}\boldsymbol{a}^H\bX)$ is linear in $\bX$~\cite{candes2013phaselift}, the AF is quartic in $\bx$, so the lifted problem is still a PR problem, now over the positive semidefinite (PSD) cone: find $\bX\succeq 0$ of rank one with $|\B(\bX)| = \bb$.

Let $\bx$ occupy at most $S$ contiguous DFT bins; because a modulation is trivial (T4), the band can be taken centered. Let $\bQ\in\C^{N\times S}$ hold the DFT columns $e^{2\pi i nk/N}/\sqrt{N}$ for $k=-S/2,\dots,S/2-1$, so that $\bx=\bQ\bc$ with $\bc\in\C^S$ and $\bX = \bQ\bW\bQ^H$ with $\bW = \bc\bc^H\in\Hb^S$. The band-limited lift $\B_Q(\bW) = \B(\bQ\bW\bQ^H)$ has $S^2$ unknowns instead of $N^2$, and $\bQ\bW\bQ^H\succeq 0$ exactly when $\bW\succeq 0$.

The trivial ambiguities T1--T4 of $\bx$ are collected in the set $\mathcal{T}(\bx)$. For almost all pulses with $S\le N/2$, $3S-1$ AF samples determine $\bx$ up to $\mathcal{T}(\bx)$~\cite[Prop.~1]{pinilla2024iaae}. We take $S=N/2$, the widest band covered by this uniqueness result. Since $\bB_{0,0}=\bI$, the origin sample fixes the trace: $\tr(\bW) = \|\bx\|_2^2 = \bb[0,0]$ for every exact fit.

\subsection{Convex Program and Recovery Guarantees}\label{sec:convex}

Let $\bW_\star = \bc\bc^H$ denote the lifted pulse in the band. The convex program takes as input an \emph{anchor} $\bW_0 = \bc_0\bc_0^H$, a rank-one estimate of $\bW_\star$ computed from the data by the initialization of Section~\ref{sec:anchor}. Given the anchor, and following~\cite{pinilla2025wavemax}, we relax the equalities $|\B_Q(\bW)|=\bb$ to inequalities and, among the matrices that satisfy them, recover $\bW_\star$ as the one best aligned with the anchor:
\begin{equation}
  \underset{\bW\in\Hb^S}{\text{maximize}}\ \ \langle \bW_0, \bW\rangle
  \ \ \text{subject to}\ \ |\B_Q(\bW)|\le \bb,\ \bW\succeq 0 .
  \tag{P}\label{eq:P}
\end{equation}
The objective is linear and the constraints are convex, so \eqref{eq:P} is convex~\cite[Lemma~1]{pinilla2025wavemax}. The lift of every member of $\mathcal{T}(\bx)$ that remains in the band is feasible with all constraints active, so the anchor is what singles out $\bW_\star$.

To state when $\bW_\star$ is the unique solution of~\eqref{eq:P}, let $\bu = \sgn(\B_Q(\bW_\star))$ hold the phases of the AF samples, and define on $\Hb^S$ the real-linear map $\mathcal{L}(\boldsymbol{Z}) = \re(\overline{\bu}\odot\B_Q(\boldsymbol{Z})) \in\R^m$, with adjoint $\mathcal{L}^*$. Since $\mathcal{L}(\bW)\le|\B_Q(\bW)|$, with equality $\mathcal{L}(\bW_\star)=\bb$ at the pulse, every feasible $\bW$ satisfies $\mathcal{L}(\bW-\bW_\star)\le\boldsymbol{0}$.

Moving the pulse from $\bc$ to $\bc+\boldsymbol{h}$ changes $\bc\bc^H$, to first order, by $\bc\boldsymbol{h}^H+\boldsymbol{h}\bc^H$; these matrices form the tangent space $E = \{\bc\boldsymbol{h}^H + \boldsymbol{h}\bc^H : \boldsymbol{h}\in\C^S\}$ of the rank-one matrices at $\bW_\star$~\cite[Sec.~IV-C]{pinilla2025wavemax}. Its orthogonal complement $E^\perp$ holds the Hermitian matrices whose range is orthogonal to $\bc$, with projection $P_{E^\perp}(\boldsymbol{Z}) = \boldsymbol{\Pi}\boldsymbol{Z}\boldsymbol{\Pi}$, $\boldsymbol{\Pi} = \bI - \bc\bc^H/\|\bc\|_2^2$. We write $\boldsymbol{Z}_E = P_E(\boldsymbol{Z})$, and likewise for the other subspaces.

The map $\mathcal{L}$ cannot be injective on all of $E$: a continuous trivial ambiguity, such as a fractional translation, moves $\bW_\star$ along a curve of rank-one matrices on which $|\B_Q|$ is constant, so the tangent of this curve lies in $K = E\cap\ker\mathcal{L}$. Injectivity is therefore required only on $E' = E\ominus K$, the orthogonal complement of $K$ in $E$:
\begin{equation}
  \alpha = \min_{\boldsymbol{Z}\in E',\ \|\boldsymbol{Z}\|_F = 1}\|\mathcal{L}(\boldsymbol{Z})\|_1 .
  \label{eq:alpha}
\end{equation}
In the opposite direction, $\|\bQ^H\bB_{\tau,f}\bQ\|_2\le1$ gives $\|\mathcal{L}(\boldsymbol{Z})\|_1\le m\|\boldsymbol{Z}\|_*$ for every $\boldsymbol{Z}\in\Hb^S$.

\begin{theorem}\label{thm:unique}
Assume (A1) $\alpha>0$; (A2) there are $\blam\in\R^m$ with $\lambda_j\ge\lambda_{\min}>0$ and $\nu,\zeta\ge0$ such that $\bY = \mathcal{L}^*(\blam) - \bW_0$ satisfies $P_{E^\perp}(\bY) \succeq \nu\boldsymbol{\Pi}$ and $\|P_{E'}(\bY)\|_F\le\zeta$; and (A3) the anchor is aligned, $P_K(\bW_0) = \boldsymbol{0}$. If $m\lambda_{\min}<\nu$ and $\zeta<\alpha\lambda_{\min}$, then $\bW_\star$ is the unique solution of~\eqref{eq:P}.
\end{theorem}

\begin{IEEEproof}
Let $\bW\neq\bW_\star$ be feasible and $\boldsymbol{H} = \bW-\bW_\star$. Feasibility gives $\mathcal{L}(\boldsymbol{H}) \le\boldsymbol{0}$ and $\boldsymbol{H}_{E^\perp} = \boldsymbol{\Pi}\bW\boldsymbol{\Pi}\succeq 0$. Since $\langle\bW_0,\boldsymbol{H}\rangle = \blam^T\mathcal{L}(\boldsymbol{H}) - \langle\bY,\boldsymbol{H}\rangle$, we bound the two terms. First, as $\mathcal{L}(\boldsymbol{H}) \le\boldsymbol{0}$ and $\mathcal{L}(\boldsymbol{H}_K) = \boldsymbol{0}$,
\begin{align*}
  \blam^T\mathcal{L}(\boldsymbol{H}) &\le -\lambda_{\min}\|\mathcal{L}(\boldsymbol{H})\|_1\\
  &\le -\lambda_{\min}\big(\alpha\|\boldsymbol{H}_{E'}\|_F - m\|\boldsymbol{H}_{E^\perp}\|_*\big).
\end{align*}
Second, $P_K(\mathcal{L}^*(\blam)) = \boldsymbol{0}$ because the range of $\mathcal{L}^*$ is orthogonal to $\ker\mathcal{L}\supseteq K$, so (A3) gives $P_K(\bY)=\boldsymbol{0}$ and, by (A2) and $\boldsymbol{H}_{E^\perp}\succeq0$,
\begin{align*}
  -\langle\bY,\boldsymbol{H}\rangle &= -\langle P_{E'}(\bY),\boldsymbol{H}_{E'}\rangle
  - \langle P_{E^\perp}(\bY),\boldsymbol{H}_{E^\perp}\rangle\\
  &\le \zeta\|\boldsymbol{H}_{E'}\|_F - \nu\|\boldsymbol{H}_{E^\perp}\|_* .
\end{align*}
Hence $\langle\bW_0,\boldsymbol{H}\rangle\le -(\alpha\lambda_{\min}-\zeta)\|\boldsymbol{H}_{E'}\|_F -(\nu-m\lambda_{\min})\|\boldsymbol{H}_{E^\perp}\|_*$, which is negative unless $\boldsymbol{H}\in K$. In that case $\boldsymbol{H} = \bc\boldsymbol{h}^H+\boldsymbol{h}\bc^H$ with $\bW_\star+\boldsymbol{H}\succeq0$; writing $\boldsymbol{h} = \beta\bc+\boldsymbol{r}$ with $\boldsymbol{r}\perp\bc$, the $2\times2$ compression of $\bW_\star+\boldsymbol{H}$ to $\{\bc,\boldsymbol{r}\}$ has determinant $-\|\bc\|_2^2\|\boldsymbol{r}\|_2^2$, so $\boldsymbol{r}=\boldsymbol{0}$ and $\boldsymbol{H} = t\,\bc\bc^H$ with $t$ real. Then $\mathcal{L}(\boldsymbol{H}) = t\,\bb = \boldsymbol{0}$ forces $t=0$, a contradiction.
\end{IEEEproof}

Conditions (A1) and (A2) are the restricted injectivity and approximate dual certificate of~\cite[Thm.~2]{pinilla2025wavemax}, stated modulo the kernel $K$ that the continuous ambiguities create; (A3) is the first-order optimality of the alignment between the anchor and the orbit of $\bW_\star$. When an exact certificate is available, the injectivity and alignment conditions are not needed:

\begin{proposition}\label{prop:cert}
If there is $\blam>\boldsymbol{0}$ such that $\bY = \mathcal{L}^*(\blam) - \bW_0$ satisfies $\bY\succeq 0$, $\bY\bc = \boldsymbol{0}$ and $\operatorname{rank}(\bY) = S-1$, then $\bW_\star$ is the unique solution of~\eqref{eq:P}.
\end{proposition}

\begin{IEEEproof}
For a feasible $\bW$, $\langle \bW_0,\bW\rangle = \blam^T\mathcal{L}(\bW) - \langle \bY,\bW\rangle\le\blam^T\bb = \langle\bW_0,\bW_\star\rangle$, using $\mathcal{L}(\bW)\le\bb$, $\langle\bY,\bW\rangle\ge0$ and $\bY\bc=\boldsymbol{0}$. Equality forces $\langle\bY,\bW\rangle=0$, so $\bW = t\,\bc\bc^H$ with $t\ge0$ because the null space of $\bY$ is spanned by $\bc$, and $\blam^T\mathcal{L}(\bW)=\blam^T\bb$ forces $t=1$.
\end{IEEEproof}

Since $\|\bv\|_1\ge\|\bv\|_2$, the smallest singular value of $\mathcal{L}$ on $E'$ bounds $\alpha$ from below; Section~\ref{sec:exp} evaluates it numerically (Table~\ref{tab:alpha}).

\section{The \textit{wavex} Algorithm}

\textit{wavex} recovers the pulse in three steps, summarized in Algorithm~\ref{alg:wavex}: an initialization that computes the anchor $\bW_0$ of~\eqref{eq:P}, convex surrogate steps that approach the solution of~\eqref{eq:P} in the lifted space, and a polish of the resulting pulse in the full signal space.

\subsection{Initialization}\label{sec:anchor}

The lifted pulse $\bW_\star$ belongs to two sets, the matrices that reproduce the measured magnitudes, $|\B_Q(\bW)|=\bb$, and the rank-one PSD matrices, and the initialization looks for a common point by moving between them. Let $P_M(\bW) = \bW + \B_Q^+\big(\bb\odot\sgn(\B_Q(\bW)) - \B_Q(\bW)\big)$ fit the measured magnitudes while keeping the current phases, with $\B_Q^+$ the minimum-norm pseudoinverse, and let $P_S(\bW) = \lambda_1\bv_1\bv_1^H$ be the best rank-one PSD approximation. Alternating them, $\bW\leftarrow P_SP_M(\bW)$, is the spectral step behind earlier estimates~\cite{pinilla2024iaae,pinilla2025wavemax}. The initialization instead uses the relaxed averaged alternating reflections (RAAR)~\cite{luke2005raar}, which reflect through both sets so that the iterate does not stall, as plain alternation can:
\begin{equation}
  \bW \leftarrow \tfrac{\beta}{2}\big(R_SR_M + \mathrm{Id}\big)(\bW) + (1-\beta)P_M(\bW),
  \quad R = 2P - \mathrm{Id},
  \label{eq:raar}
\end{equation}
with $\beta = 0.9$ and $300$ steps, from a single random start $\boldsymbol{g}\boldsymbol{g}^H$ scaled to $\|\boldsymbol{g}\|_2 = \boldsymbol{A}[0,0]^{1/4} = \|\bx\|_2$. Applying $\B_Q^+$ by ten conjugate-gradient steps from zero keeps the correction in the range of $\B_Q^*$, so the pulse remains a fixed point even where $\B_Q$ has a null space. The anchor is $\bc_0 = \lead(P_M(\bW))$.

\subsection{Convex Surrogate Steps}\label{sec:lifted-stage}

Problem~\eqref{eq:P} has $S^2 = 4096$ real unknowns and up to $16384$ second-order cone constraints, which is costly for interior-point solvers. \textit{wavex} instead applies the fast adaptive shrinkage/thresholding algorithm (FASTA)~\cite{goldstein2014fasta}, started at $\bW_0$, to the penalized form
\begin{equation}
  \underset{\bW\succeq 0}{\text{minimize}}\ \ \tfrac12\big\|\,|\B_Q(\bW)| + \boldsymbol{\xi} - \bb\,\big\|_2^2 + \mu\tr(\bW) .
  \label{eq:pen}
\end{equation}
The fit penalizes violations of the constraints of~\eqref{eq:P}, and the slack $\boldsymbol{\xi}\ge\boldsymbol{0}$, refreshed between rounds as $\boldsymbol{\xi}=\max(\bb-|\B_Q(\bW)|,\boldsymbol{0})$, stops pulling up the samples that already lie below their bound. The trace, which is the nuclear norm on the PSD cone, is a convex surrogate of rank; since $\tr(\bW) = \bb[0,0]$ on exact fits, it does not change the set of exact solutions, and its $O(\mu)$ bias is removed by continuation, $\mu\leftarrow0.1\mu$ over three rounds, from $\mu = 10^{-2}\|\B_Q^*(\bb)\|_*$.

With $\bt=\bb-\boldsymbol{\xi}$ and $\bu_k = \sgn(\B_Q(\bW_k))$ at the iterate $\bW_k$, the convex quadratic $G_k(\bW) = \tfrac12\|\B_Q(\bW) - \bt\odot \bu_k\|_2^2$ majorizes the fit, because $|z - tu| \ge \big||z|-t\big|$ for $|u|=1$, and touches it with the same gradient at $\bW_k$. The forward--backward step of FASTA at $\bW_k$ is therefore a step on the convex semidefinite program of minimizing $G_k(\bW)+\mu\tr(\bW)$ over $\bW\succeq0$, a majorization--minimization scheme, whose proximal map soft-thresholds the eigenvalues. FASTA chooses the step $\gamma$ adaptively~\cite{goldstein2014fasta}, and each iteration costs one $S\times S$ Hermitian eigendecomposition and two applications of $\B_Q$. The lifted estimate is $\bx_{\rm L} = \bQ\,\lead(\bW)$.

\subsection{Polish}\label{sec:polish}

The lifted estimate need not be the pulse itself: every member of $\mathcal{T}(\bx)$ is an equally valid solution, and a member modulated by a fraction of a bin (T4) lies outside the band, where the convex steps can only return its band-limited approximation. The polish (i.e. step to remove the trivial ambiguity) moves the estimate onto that member in the full space $\C^N$ by decreasing the intensity misfit
\begin{equation}
  \ell(\bq) = \tfrac12\sum_{(\tau,f)\in\Omega}\big(|\bq^H\bB_{\tau,f}\,\bq|^2 - \boldsymbol{A}[\tau,f]\big)^2,
  \quad \bq\in\C^N,
  \label{eq:polish}
\end{equation}
with at most 150 Levenberg--Marquardt (LM) iterations~\cite{marquardt1963} started at $\bq=\bx_{\rm L}$; the result is the estimate $\hat\bx$.

\begin{algorithm}[t]
\caption{\textit{wavex}: estimate $\hat\bx$ from AF samples $\boldsymbol{A}$ on $\Omega$}\label{alg:wavex}
\begin{algorithmic}[1]
\Require $\boldsymbol{A}$ on $\Omega$; $S=N/2$; $\beta=0.9$
\State $\bb\gets\sqrt{\boldsymbol{A}}$; $\bW\gets\boldsymbol{g}\boldsymbol{g}^H$, $\boldsymbol{g}\in\C^S$ random, $\|\boldsymbol{g}\|_2=\boldsymbol{A}[0,0]^{1/4}$
\Statex \textit{Initialization}
\For{$300$ steps}
  \State $\bW_M\gets P_M(\bW)$ \Comment{fit magnitudes, keep phases}
  \State $\bW\gets\tfrac{\beta}{2}\big(R_SR_M+\mathrm{Id}\big)(\bW)+(1-\beta)\bW_M$
\EndFor
\State $\bc_0\gets\lead(P_M(\bW))$; $\bW\gets\bc_0\bc_0^H$ \Comment{anchor $\bW_0$}
\State $\boldsymbol{\xi}\gets\boldsymbol{0}$; $\mu\gets10^{-2}\|\B_Q^*(\bb)\|_*$
\Statex \textit{Convex surrogate steps}
\For{$3$ rounds}
  \Repeat \Comment{FASTA with step $\gamma$}
    \State $\bz\gets\B_Q(\bW)$
    \State $\bW\gets\bW-\gamma\,\B_Q^*\big(\sgn(\bz)\odot(|\bz|+\boldsymbol{\xi}-\bb)\big)$
    \State $\bV\bD\bV^H\gets\bW$ \Comment{eigendecomposition}
    \State $\bW\gets\bV\max(\bD-\gamma\mu\bI,\boldsymbol{0})\bV^H$
  \Until{the residual falls below $10^{-5}$ of its maximum}
  \State $\boldsymbol{\xi}\gets\max(\bb-|\B_Q(\bW)|,\boldsymbol{0})$; $\mu\gets0.1\mu$
\EndFor
\State $\bx_{\rm L}\gets\bQ\,\lead(\bW)$ \Comment{lifted estimate}
\Statex \textit{Polish}
\State $\hat\bx\gets$ LM iterations on $\ell$ in~\eqref{eq:polish}, started at $\bx_{\rm L}$
\State \Return $\hat\bx$ (noiseless data) or $\bx_{\rm L}$ (noisy data)
\end{algorithmic}
\end{algorithm}

\section{Numerical Experiments}\label{sec:exp}

Following~\cite{pinilla2024iaae}, we build 100 unit-norm pulses with $N=128$ samples on $t\in[-1,1]$~$\mu$s ($\Delta t = 15.7$~ns), each as the inverse DFT of a Gaussian-shaped amplitude spectrum multiplied by a uniformly random phase. The centered band of $N/2$ bins holds all but $1.1\times10^{-7}$ of the energy of each pulse, so the pulses are band-limited with $S\le N/2$ up to the Gaussian tail. The error is $\operatorname{dist}(\bx,\hat\bx)/\|\bx\|_2$, with the distance of~\cite[Def.~2]{pinilla2024iaae} taken over $\mathcal{T}(\bx)$ with continuous shifts and modulations, and a trial succeeds when it is below $10^{-3}$. \textit{wavex} uses at most $1000$ FASTA iterations per round and a single random start per pulse; all results use one seed.

\begin{table}[h]
\centering
\caption{Numerical validation of (A1) on pulses 1, 6 and 11: singular values of $\mathcal{L}$ on $E$ (range over the pulses)}
\label{tab:alpha}
\setlength{\tabcolsep}{5pt}
\begin{tabular}{@{}rrccc@{}}
\toprule
Delays kept & $m$ & On $K$ (shift) & Smallest on $E'$ & Next on $E'$ \\
\midrule
All        & 16384 & $\le3.7\!\times\!10^{-13}$ & $5.2$--$10.7\!\times\!10^{-3}$ & $3.1$--$3.6$ \\
Every 2nd  & 8192  & $\le2.5\!\times\!10^{-13}$ & $3.7$--$7.6\!\times\!10^{-3}$  & $2.2$--$2.5$ \\
Every 3rd  & 5504  & $\le2.1\!\times\!10^{-13}$ & $3.0$--$6.2\!\times\!10^{-3}$  & $1.8$--$2.1$ \\
Every 5th  & 3328  & $\le1.6\!\times\!10^{-13}$ & $2.3$--$4.8\!\times\!10^{-3}$  & $1.4$--$1.6$ \\
Every 7th  & 2432  & $\le1.5\!\times\!10^{-13}$ & $2.0$--$4.0\!\times\!10^{-3}$  & $1.2$--$1.3$ \\
Every 8th  & 2048  & $\le1.3\!\times\!10^{-13}$ & $1.6$--$2.7\!\times\!10^{-3}$  & $0.13$--$0.16$ \\
\bottomrule
\end{tabular}
\end{table}

\emph{Numerical validation of (A1):} Table~\ref{tab:alpha} evaluates condition (A1) of Theorem~\ref{thm:unique} on three of the pulses, for the complete AF and the delay decimations of Table~\ref{tab:incomplete}, through the singular values of $\mathcal{L}$ on $E$; the smallest on $E'$ is a lower bound for $\alpha$. In every case $\alpha>0$: $K$ is one-dimensional (the fractional shift), and $E'$ has one weak direction, the modulation by a fraction of a bin that the band cannot express and the polish corrects, with all other singular values above one up to every seventh delay. Keeping every eighth delay adds two more weak directions, consistent with its failure in Table~\ref{tab:incomplete}.

\begin{figure}[h]
  \centering
  \includegraphics[width=\columnwidth]{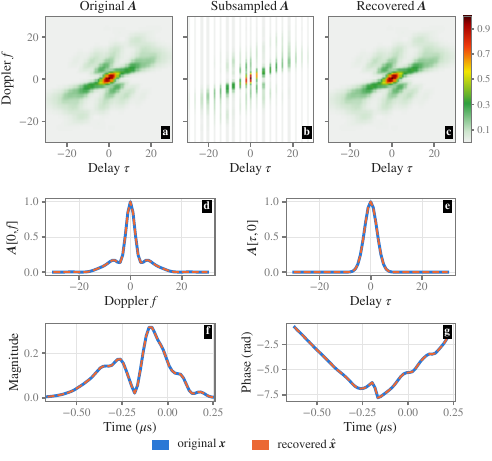}
  \caption{Recovery by \textit{wavex} when every third delay is kept (66\% of the AF samples removed). (a)--(c) Original, subsampled and recovered AFs, normalized to a unit peak. (d), (e) Zero-delay and zero-Doppler slices. (f), (g) Magnitude and phase of $\hat\bx$ (red) and $\bx$ (blue), aligned over $\mathcal{T}(\bx)$.}
  \label{fig:example}
\end{figure}

\emph{Example:} Fig.~\ref{fig:example} shows a recovery by \textit{wavex} when only every third delay is measured. The AF of the recovered pulse $\hat\bx$ matches the original over the whole delay-Doppler plane, including the unmeasured delays, and $\hat\bx$ overlaps $\bx$ after alignment over $\mathcal{T}(\bx)$, with a relative error of $4.8\times10^{-6}$ ($9.8\times10^{-4}$ for the anchor).

\begin{figure}[t]
	\centering
	\includegraphics[width=\columnwidth]{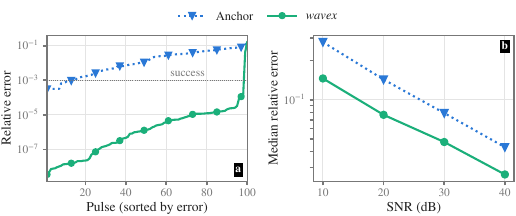}
	\caption{Relative error of the anchor and of \textit{wavex}. (a) Sorted errors of the 100 pulses from the complete AF; the dotted line marks the success threshold $10^{-3}$. (b) Median error over 20 pulses versus SNR on the complete AF.}
	\label{fig:results}
\end{figure}

\emph{Complete AF:} Fig.~\ref{fig:results}(a) shows that the anchor has a median error of $1.4\times10^{-2}$ and that \textit{wavex} recovers 98 of the 100 pulses, with a median error of $1.4\times10^{-6}$. A recovery can be verified from the data alone, since the final misfit~\eqref{eq:polish} of a recovered pulse is about $10^{-8}$. The median run time per pulse on one CPU thread is 495~s.

\emph{Incomplete AF:} Table~\ref{tab:incomplete} uses 20 pulses (pulses $0,5,\dots,95$) and removes delays or Doppler bins, as in Fig.~\ref{fig:example}. When every $L$-th delay is kept with $L$ coprime to $N$, or the central half of the delays or of the Doppler bins is kept ($|\tau|<N/4$ or $|f|<N/4$), \textit{wavex} recovers every pulse. When $L$ divides $N$ ($L=2,4,8$ for $N=128$), the success rate drops to 55\%, 25\% and 0\%. Keeping only the middle half of the delays or of the Doppler bins removes the origin rows, which the uniqueness proof of~\cite{pinilla2024iaae} requires, and recovery fails.

\begin{table}[h]
\centering
\caption{Median error and Success Rate of \textit{wavex} for Incomplete AF}
\label{tab:incomplete}
\setlength{\tabcolsep}{6pt}
\begin{tabular}{@{}lrcc@{}}
\toprule
Measured & Removed & Median error & Success \\
\midrule
All delays             & 0\%     & $6.2\!\times\!10^{-7}$ & 100\% \\
Every 3rd delay        & 66\%    & $1.2\!\times\!10^{-7}$ & 100\% \\
Every 5th delay        & 80\%    & $4.5\!\times\!10^{-7}$ & 100\% \\
Every 7th delay        & 86\%    & $7.8\!\times\!10^{-7}$ & 100\% \\
Every 2nd delay        & 50\%    & $3.6\!\times\!10^{-4}$ & 55\%  \\
Every 4th delay        & 75\%    & $1.7\!\times\!10^{-3}$ & 25\%  \\
Every 8th delay        & 88\%    & $5.9\!\times\!10^{-2}$ & 0\%   \\
$|\tau| < N/4$         & 50\%    & $1.4\!\times\!10^{-7}$ & 100\% \\
$|f| < N/4$            & 50\%    & $6.5\!\times\!10^{-7}$ & 100\% \\
Middle half of delays  & 50\%    & $2.6\!\times\!10^{-1}$ & 10\%  \\
Middle half of Doppler & 50\%    & $7.4\!\times\!10^{-1}$ & 5\%   \\
\bottomrule
\end{tabular}
\end{table}

\emph{Noisy AF:} We add Gaussian noise $\be$ to the complete AF with $\text{SNR}=10\log_{10}(\|\boldsymbol{A}\|_F^2/\|\be\|_2^2)$. For noisy data \textit{wavex} returns $\bx_{\rm L}$, because the polish fits the noise. Fig.~\ref{fig:results}(b) shows that its error degrades gracefully, with medians of 0.146, 0.077, 0.047 and 0.027 at 10, 20, 30 and 40~dB, about half those of the anchor (0.279, 0.144, 0.079 and 0.043).

\section{Conclusion}

This letter showed that radar waveform design from the magnitude of the conventional AF, so far addressed by non-convex solvers without global guarantees, admits a convex formulation. Lifting the problem to PSD matrices on the pulse band makes each AF sample the magnitude of a linear functional, and recovery becomes a convex program with a dual certificate as a sufficient condition for exactness. The proposed \textit{wavex} combines an initialization by relaxed reflections, convex surrogate steps and a local polish. It accurately recovered band-limited pulses from complete AFs and from AFs with most samples removed, provided the remaining ones determine the pulse, and degraded gracefully with noise. 

\newpage
\bibliographystyle{IEEEtran}
\bibliography{refs}

\end{document}